\documentclass{amsart}
\usepackage{amsmath, amssymb, amsthm, amsfonts, amsxtra}
\usepackage{enumitem}
\usepackage{graphicx, float, caption, subcaption} 
\usepackage{array,longtable} 
\usepackage{verbatim}
\usepackage{cancel}
\usepackage[dvipsnames]{xcolor}
\usepackage{ifthen, mathtools, hhline, bbm} 
\usepackage[mode=image|tex]{standalone} 
\usepackage{float}
\usepackage{tikz}
\usetikzlibrary{arrows.meta,backgrounds,shapes,calc} 
\usepackage[bookmarks=true, bookmarksopen=false,%
    colorlinks=true,%
    linkcolor=darkblue,%
    citecolor=darkblue,%
    filecolor=darkblue,%
    menucolor=darkblue,%
    linktoc=all,%
    urlcolor=darkblue
]{hyperref}
\usepackage[capitalize,noabbrev]{cleveref}

\definecolor{darkblue}{cmyk}{1,0.4,0,0.4}  

\newtheorem{proposition}{Proposition}[section]
\newtheorem{theorem}[proposition]{Theorem}

\newtheorem{lemma}[proposition]{Lemma}

\theoremstyle{definition}

\theoremstyle{remark}

\numberwithin{equation}{section}

\crefname{proposition}{Proposition}{Propositions}
\crefname{theorem}{Theorem}{Theorems}
\crefname{corollary}{Corollary}{Corollaries}
\crefname{lemma}{Lemma}{Lemmas}
\crefname{definition}{Definition}{Definitions}
\crefname{remark}{Remark}{Remarks}
\crefname{table}{Table}{Tables}

\Crefname{proposition}{Proposition}{Propositions}
\Crefname{theorem}{Theorem}{Theorems}
\Crefname{corollary}{Corollary}{Corollaries}
\Crefname{lemma}{Lemma}{Lemmas}
\Crefname{definition}{Definition}{Definitions}
\Crefname{remark}{Remark}{Remarks}
\crefname{table}{Table}{Tables}

\newcommand{\Kekule}{\ensuremath{\mathfrak{K}}}

\begin{document}

\title{Alternative Approaches for Counting Weakly Increasing Matrices}
\author{Leo Yicheng Yang}
\address{North Carolina School of Science and Mathematics, Durham, North Carolina, USA}
\email{leo.yc.yang@gmail.com}
\date{July 1, 2025}

\maketitle

\begin{abstract}
This paper presents two alternative approaches for counting the number of two-row weakly increasing matrices,
which are $2\times n$ matrices whose entries are integers from $1$ to $k$ and are weakly increasing along all rows and columns, for any positive integers $n$ and $k$.
The first approach establishes a bijection between the set of such matrices and the set of Kekul\'e structures for certain hexagonal benzenoids.
The second approach reduces the problem to counting the number of pairs of non-intersecting lattice paths.
These approaches reveal interesting connections between combinatorial problems that arise in different domains.
\end{abstract}

\section{Problem Introduction}
\label{sec:intro}

An $m \times n$ matrix $A = \{a_{i,j}\}$ is \emph{weakly increasing}
if its entries increase weakly along each row and down each column,
i.e., $a_{i,j} \le a_{i',j'}$ for any $1 \le i \le i' \le m$ and $1 \le j \le j' \le n$.
Given natural numbers $n$ and $k$, let $\mathfrak{M}_{n, k}$ denote the set of $2 \times n$ weakly increasing matrices with entries in $\{1, \ldots, k\}$.
We are interested in $|\mathfrak{M}_{n, k}|$ as a function of $n$ and $k$,
which was posed as Problem \#12497 in the \textit{American Mathematical Monthly}.

For example, there are six matrices in $\mathfrak{M}_{2, 2}$:
$$
\left[
\begin{array}{cc}
1 & 1 \\
1 & 1
\end{array}
\right]
\quad
\left[
\begin{array}{cc}
1 & 1 \\
1 & 2
\end{array}
\right]
\quad
\left[
\begin{array}{cc}
1 & 1 \\
2 & 2
\end{array}
\right]
\quad
\left[
\begin{array}{cc}
1 & 2 \\
1 & 2
\end{array}
\right]
\quad
\left[
\begin{array}{cc}
1 & 2 \\
2 & 2
\end{array}
\right]
\quad
\left[
\begin{array}{cc}
2 & 2 \\
2 & 2
\end{array}
\right]
$$
As another example, here is a matrix in $\mathfrak{M}_{6,7}$:
\begin{align}\label{eq:matrix-example}
W =
\begin{bmatrix}
1 & 1 & 2 & 3 & 6 & 6\\
1 & 1 & 2 & 4 & 6 & 7
\end{bmatrix}
\end{align}

From Chen et al.~\cite{ChenRamkumarYang25}, the formula for $|\mathfrak{M}_{n, k}|$ is:
\begin{align}\label{eq:count-formula}
|\mathfrak{M}_{n, k}| = \frac{{n+k-1 \choose k-1}{n+k \choose k-1}}{k}
\end{align}
The authors derived this formula from Macmahon's formula for counting the number of plane partitions~\cite{MacMahon1916Combinatory} within a $2\times n \times k$ grid,
by establishing a bijection between $\mathfrak{M}_{n, k}$ and the set of plane partitions.

This paper presents two alternative approaches for deriving $|\mathfrak{M}_{n, k}|$.
To the best of our knowledge, these approaches are novel.
We believe that they reveal interesting connections between combinatorial problems that arise in different domains.

The first approach (\cref{sec:benzenoid}) draws a bijection between $\mathfrak{M}_{n, k}$ and the set of Kekul\'e structures of hexagonal benzenoids,
whose count has been established by Cyvin in 1986~\cite{Cyvin1986} (see also~\cite{cyvin1988kekulé}).
A hexagonal benzenoid is made up of regular hexagonal benzene rings arranged in a hexagonal structure,
and a Kekul\'e structure is a selection of edges in this structure such that every vertex is paired with exactly one other vertex.
We observe that the positions of vertically oriented edges among the selected edges follow certain patterns and uniquely identify each Kekul\'e structure.
Furthermore, there is a bijection between these sets of positions and $\mathfrak{M}_{n, k}$.

The second approach (\cref{sec:paths}) draws a bijection between $\mathfrak{M}_{n, k}$ and the set of pairs of non-intersecting paths in a two-dimensional rectangular lattice,
where edges only point upward and rightward.
The locations of upward edges in each path determine which entries in a row see strict increases and by how much.
The Lindstr\"{o}m-Gessel-Viennot Lemma, which counts the number of non-intersecting paths in a directed graph, can be applied to produce the formula in~\Cref{eq:count-formula}.

Finally, in \Cref{sec:conclude}, we comment on how these approaches can be naturally extended to solve the more general problem
of counting the number of weakly increasing matrices of dimension $m \times n$.

\section{Connection to Kekul\'{e} Structures for Hexagonal Benzenoids}
\label{sec:benzenoid}

\begin{figure}[t]
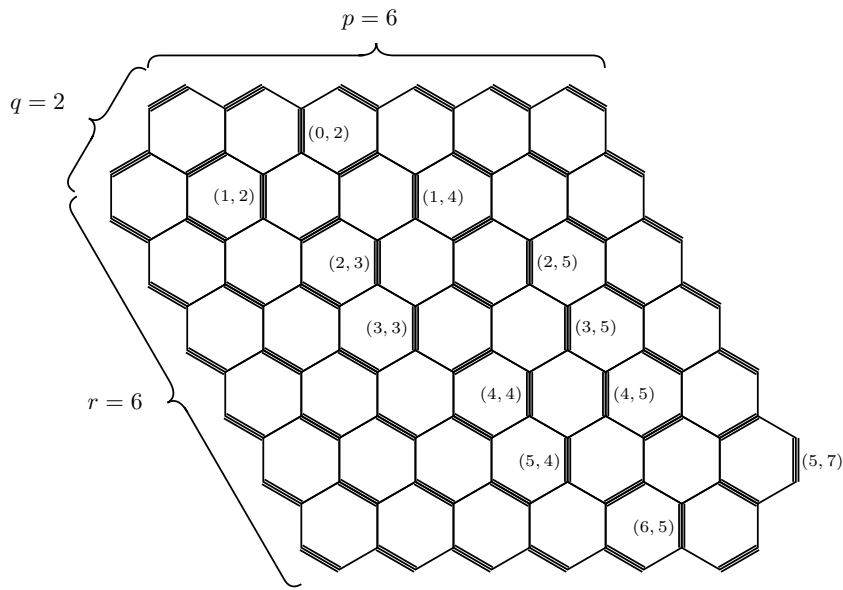

    \centering
    \includestandalone[scale=0.9]{figures/rawkekule}
    \caption{\label{fig:kekule}Hexagonal benzenoid $O\{p,q,r\}$ with dimensions $p=6$, $q=2$, and $r=6$,
    with a particular Kekul\'{e} structure whose selected edges are shown as triple lines.
    The v-bars (selected edges that are vertical) in the Kekul\'{e} structure are further marked with their positions.}
\end{figure}

A \emph{hexagonal benzenoid} $O\{p, q, r\}$ consists of $(p+1)(q+r-1)-2$ regular hexagonal benzene rings (units) arranged in a hexagonal structure, as illustrated in \Cref{fig:kekule}.
A Kekul\'{e} structure is a perfect matching of $O$, i.e.,
a subset of edges such that every vertex is paired with exactly one other vertex.
We will refer to edges in this subset as \emph{selected} edges.

Each selected edge in the Kekul\'{e} structure represents a double bond between two carbon atoms in a benzene ring (C6H6).
Let $\Kekule_{p,q,r}$ denote the set of Kekul\'{e} structures for a hexagonal benzenoid with dimensions $(p, q, r)$.
An explicit combinatorial formula for $| \Kekule_{p,q,r} |$ already exists,
having been derived by Cyvin in 1986~\cite{Cyvin1986} (see also~\cite{cyvin1988kekulé}):
\begin{align}\label{eq:kekule}
| \Kekule_{p,q,r} | = \prod_{i=0}^{q-1} \frac{{p+r+i \choose r}}{{r+i \choose r}}
\end{align}
Notice that plugging in $q=2$ yields:
\begin{align*}
| \Kekule_{p,2,r} |
= \frac{{p+r \choose r} {p+r+1 \choose r}}{{r \choose r}{r+1 \choose r}}
= \frac{{p+r \choose r} {p+r+1 \choose r}}{r+1}
\end{align*}
which is precisely the formula for $| \mathfrak{M}_{n,k} |$ in \Cref{eq:count-formula} after substituting $r+1$ for $k$ and $p$ for $n$. 
Hence, to show \Cref{eq:count-formula}, it suffices to establish the following:

\begin{theorem}\label{thm:benzenoid}
  There exists a bijection between $\mathfrak{M}_{n, k}$ and $\Kekule_{n,2,k-1}$.
\end{theorem}

We outline the proof in the remainder of this section.
For detailed proofs of all theorems and lemmas, see \Cref{app:proof}.

\subsection{Analysis of Kekul\'{e} Structures}
\label{sec:benzenoid:kekule}

Unless otherwise noted, $K$ in this subsection refers to any Kekul\'{e} structure in $\Kekule_{n,2,k-1}$. Define a \emph{v-bar} as a selected vertical edge in $K$.
To denote the location of a vertical edge (and potential v-bar) in $K$, we use $(i,j)$,
where $i \in [0,r]$ identifies the row containing the vertical edge (beginning with $0$ for the top row of units),
and $j \in [0,n+1]$ is the number of vertical edges in the same row strictly to its left.
Note that for $i=0$ and $i=r$ (top and bottom rows), $j \le n$, because these two rows contain $n$ units each instead of $n+1$ units.

A key insight is that the set of v-bars uniquely identifies a Kekul\'{e} structure, and must follow certain constraints.
To aid in our analysis, we define several useful substructures of $K$.

\begin{figure}[t]
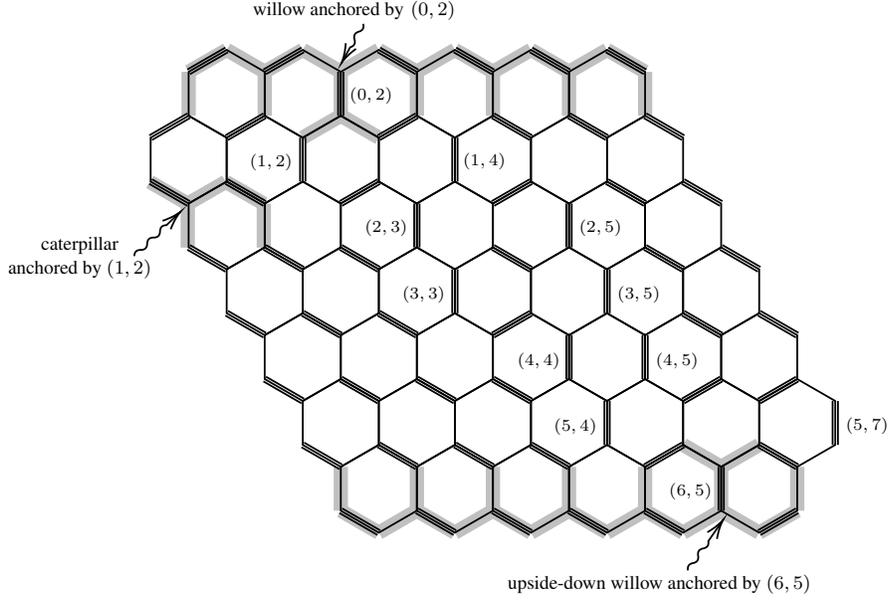

  \centering
  \includestandalone[scale=0.9]{figures/filledkekule}
  \caption{\label{fig:substructures}Substructures of $K$, willows and caterpillars, highlighted in gray. ``Leftover" edges outlined in \cref{lemma:bottom-row} are not colored.}
\end{figure}

A \emph{willow}, anchored by v-bar $(i,j)$ and spanning $[j_1,j_2]$ (where $j \in [j_1,j_2]$), is an assignment of an edge subset to either ``selected'' or ``unselected'' defined as follows, illustrated in \Cref{fig:substructures}:
\begin{itemize}
\item It contains all vertical edges in row $i$ from $(i,j_1)$ to $(i,j_2)$, inclusive;
  all are unselected except the anchor $(i,j)$, which is selected.
\item It contains the shortest path from the top of $(i,j)$ to the top of $(i,j_1)$,
  and the shortest path from the top of $(i,j)$ to the top of $(i,j_2)$; the edges of each path alternate between unselected and selected,
  starting from the edges incident to $(i,j)$, which are unselected.
\item It contains the two edges incident to the bottom of $(i,j)$; both are unselected. 
\end{itemize}

A \emph{caterpillar}, anchored by v-bar $(i-1,j)$ where $i \in [1,r]$, is defined as follows, illustrated in \Cref{fig:substructures}:
\begin{itemize}
\item It is anchored by a v-bar at $(i-1,j)$.
\item It contains all vertical edges in row $i$ from $(i,0)$ to $(i,j-1)$, inclusive;
  all are unselected.
\item It contains the shortest path from the bottom of $(i-1,0)$ to the top of $(i,j-1)$;
  edges in each path start as selected and alternate between unselected and selected.
\end{itemize}

We say that $K$ \emph{contains} a willow or caterpillar if $K$ agrees with the substructure in the selection of edges.

Starting from the top, we show that $K$ must contain exactly one v-bar in its top row,
exactly two v-bars in each of the $r-1$ middle rows,
and exactly one v-bar in the bottom row.
Furthermore, these v-bars must be stacked in a particular manner.
These constraints and the arguments supporting them are captured in the following sequence of lemmas. 

\begin{lemma}[Top Row]\label{lemma:top-row}
  There is exactly one v-bar in the top row of $K$,
  and $K$ contains a willow anchored at this v-bar and spanning $[0,n]$.
\end{lemma}

\begin{lemma}[One to the Right]\label{lemma:one-to-right}
  Suppose $(i-1,j')$, where $i \in [2,r-1]$, is the rightmost v-bar in row $i-1$ of $K$.
  There exists exactly one v-bar $(i,j)$ in row $i$ of $K$ such that $j \ge j'$
  (i.e., to the right of $(i-1,j')$ vertically),
  and $K$ contains a willow anchored at $(i,j)$ and spanning $[j',n+1]$.
\end{lemma}

\begin{lemma}[Row 1]\label{lemma:row-1}
  Suppose the v-bar in the top row of $K$ is at $(0,j)$ and $r>1$.
  There are exactly two v-bars $(1,j_1)$ and $(1,j_2)$ in row $1$ of $K$ such that $j_1 \le j < j_2$
  (i.e., on the two sides of $(0,j)$ vertically);
  $K$ contains a willow anchored at $(1,j_1)$ spanning $[0,j]$,
  and a willow anchored at $(1,j_2)$ spanning $[j+1,n+1]$.
\end{lemma}

\begin{lemma}[None to the Left]\label{lemma:nothing-on-left}
  Suppose $(i-1,j')$, where $i \in [2,r]$, is the leftmost v-bar in row $i-1$ of $K$.
  There exists no v-bar $(i,j)$ in row $i$ of $K$ such that $j < j'$
  (i.e., to the left of $(i-1,j')$ vertically),
  and $K$ contains a caterpillar anchored at $(i-1,j')$.
\end{lemma}

\begin{lemma}[Middle Rows]\label{lemma:middle-rows}
  If row $i-1$ of $K$, where $i \in [2,r-1]$, contains exactly two v-bars $(i-1,j'_1)$ and $(i-1,j'_2)$,
  then there are exactly two v-bars $(i,j_1)$ and $(i,j_2)$ in row $i$, with $j'_1 \le j_1 < j'_2 \le j_2$.
  Furthermore, $K$ contains a caterpillar anchored at $(i-1,j')$,
  a willow anchored at $(i,j'_1)$ spanning $[j'_1,j'_2-1]$, and a willow anchored at $(i,j_2)$ spanning $[j'_2,n+1]$.
\end{lemma}

\begin{lemma}[Bottom Row]\label{lemma:bottom-row}
  If row $r-1$ of $K$ has two v-bars $(r-1,j'_1)$ and $(r-1,j'_2)$, the bottom row of $K$ has exactly $1$ v-bar $(r,j)$ such that $j'_1 \le j < j'_2$. $K$ contains an upside-down willow anchored by this v-bar spanning $[0,n]$. The remaining edges, the shortest path from the bottom of $(r-1,0)$ to the bottom of $(r-1,n+1)$, must be selected/unselected in a unique way.
\end{lemma}

\begin{lemma}[]\label{lemma:benzenoid}
  For every $K \in \Kekule_{n,2,r}$, there are exactly $2r$ v-bars positioned at
  $(0,x_1), (1,y_1), (1, x_2+1), \ldots, (r-1, y_{r-1}), (r-1, x_r+1), (r, y_r)$,
  where $0 \le x_1 \le x_2 \le \cdots \le x_r \le n$,
  $0 \le y_1 \le y_2 \le \cdots \le y_r \le n$,
  $x_i \ge y_i$ for all $i$,
  and $y_i \le x_{i+1}$ for all $i \in [1,r-1]$.
\end{lemma}
  
\subsection{Decomposition of Weakly Increasing Matrices into Pulse Matrices}
\label{sec:benzenoid:pulse-matrices}

A $2 \times n$ \emph{pulse matrix} is a binary matrix uniquely identified by a pair $(x, y)$
where $0 \le y \le x \le n$,
such that the entries in the top row $\ell_{1,j} = 0$ for all $j \le x$ and $1$ for all $j > x$,
and the entries in the bottom row $\ell_{2,j} = 0$ for all $j \le y$ and $1$ for all $j > y$.
In other words, $x$ and $y$ represent the number of leading $0$'s in the rows.

\begin{lemma}\label{lemma:pulse}
  Each $M \in \mathfrak{M}_{n,k}$ can be written uniquely as $J_{2,n} + \sum_{i=1}^{k-1} L_i$,
  where $J_{2,n}$ is a $2 \times n$ matrix of all $1$'s,
  each $L_i$ is a pulse matrix identified by $(x_i,y_i)$,
  and $L_{i-1} \ge L_i$ for all $i \in [2,k-1]$ (where $\ge$ denotes conjunctive component-wise $\ge$).
\end{lemma}

For example, the matrix $W \in \mathfrak{M}_{6,7}$ in \Cref{eq:matrix-example} can be written as follows:
\begin{align*}
W =
\begin{bmatrix}
1 & 1 & 2 & 3 & 6 & 6\\
1 & 1 & 2 & 4 & 6 & 7
\end{bmatrix}
=&
\begin{bmatrix}
1 & 1 & 1 & 1 & 1 & 1\\
1 & 1 & 1 & 1 & 1 & 1
\end{bmatrix}\ldots J_{2,6}\\
+&
\begin{bmatrix}
0 & 0 & 1 & 1 & 1 & 1\\
0 & 0 & 1 & 1 & 1 & 1
\end{bmatrix}\ldots L_1\,\text{identified by}\,(2,2)\\
+&
\begin{bmatrix}
0 & 0 & 0 & 1 & 1 & 1\\
0 & 0 & 0 & 1 & 1 & 1
\end{bmatrix}\ldots L_2\,\text{identified by}\,(3,3)\\
+&
\begin{bmatrix}
0 & 0 & 0 & 0 & 1 & 1\\
0 & 0 & 0 & 1 & 1 & 1
\end{bmatrix}\ldots L_3\,\text{identified by}\,(4,3)\\
+&
\begin{bmatrix}
0 & 0 & 0 & 0 & 1 & 1\\
0 & 0 & 0 & 0 & 1 & 1
\end{bmatrix}\ldots L_4\,\text{identified by}\,(4,4)\\
+&
\begin{bmatrix}
0 & 0 & 0 & 0 & 1 & 1\\
0 & 0 & 0 & 0 & 1 & 1
\end{bmatrix}\ldots L_5\,\text{identified by}\,(4,4)\\
+&
\begin{bmatrix}
0 & 0 & 0 & 0 & 0 & 0\\
0 & 0 & 0 & 0 & 0 & 1
\end{bmatrix}\ldots L_6\,\text{identified by}\,(6,5)
\end{align*}

Combining \Cref{lemma:pulse} and \Cref{lemma:benzenoid}, we see a clear correspondence between between
the $x_i$'s and $y_i$'s in \Cref{lemma:benzenoid} (v-bar positions in a Kekul\'{e} structure in $\Kekule_{n,2,r}$) and
the $x_i$'s and $y_i$'s in \Cref{lemma:pulse} (number of leading $0$'s in the pulse matrices for a matrix in $\mathfrak{M}_{n,k}$),
by setting $r = k-1$.
This leads to \Cref{thm:benzenoid} (again, see \Cref{app:proof} for complete proofs).

For example, the Kekul\'{e} structure shown in \Cref{fig:kekule} maps to the weakly increasing matrix in \Cref{eq:matrix-example} under this bijection.
\section{Connection to Non-Intersecting Lattice Paths}
\label{sec:paths}

\subsection{Non-Intersecting Lattice Paths}
\label{sec:paths:lattice}

Consider a two-dimensional lattice whose vertices are points on the Cartesian plane with integer coordinates,
and edges have the forms $(x,y) \to (x+1,y)$ and $(x,y) \to (x,y+1)$ for all $(x,y)$,
i.e., edges only point upward or rightward.
We show that there is a bijection between $\mathfrak{M}_{n, k}$ and
the set of pairs of non-intersecting paths, one connecting $a_1$ to $b_1$ and the other connecting $a_2$ to $b_2$,
for strategically chosen vertices $a_1$, $a_2$, $b_1$, and $b_2$.

\begin{figure}[t]
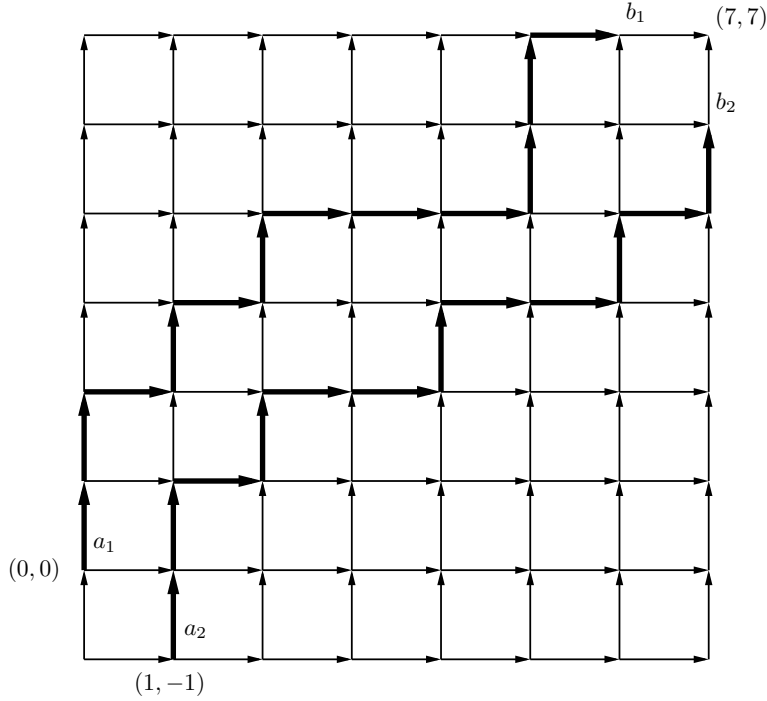

  \centering
  \includestandalone[scale=.9]{figures/2tuple}
  \caption{\label{fig:2tuple}A pair of non-intersecting lattice paths, shown in thicker lines.
  In this example, $n=6$ and $k=7$, and
  the two paths correspond to the row vectors of the example matrix in \Cref{eq:matrix-example}.}
\end{figure}

Specifically, we set $a_1 = (0,0)$, $a_2 = (1,-1)$, $b_1 = (k-1,n)$, and $b_2 = (k,n-1)$.
One example pair of non-intersecting paths is illustrated in \Cref{fig:2tuple}.
Since all edges in the lattice point upward or rightward,
and corresponding source and destination vertices always differ by $(k-1, n)$ in their location,
each path must contain exactly $k-1$ rightward edges and $n$ upward edges.

Intuitively, each path corresponds to a $1 \times n$ weakly increasing row vector $V = [v_1, \ldots, v_n]$ with entries in $[1,k]$.
To construct $V$, consider the $n+1$ path segments separated by the $n$ upward edges in the path;
each segment consists of zero or more consecutive rightward edges.
Let $\delta_i \ge 0$ ($i \in [0,n]$) denote the length of each of these rightward-only segments.
The $i$-th column of $V$ is given by the sum $v_i = 1 + \sum_{j\in[0,i)} \delta_j$.
In other words, starting from the base value of $1$, each $\delta_{i-1}$ specifies the amount of additional increase needed to obtain $v_i$.
Since $\sum_{i \in [0,n]} \delta_i = k-1$ by construction (every path has exactly $k-1$ rightward edges),
it is easy to see that $1 \le v_1 \le v_2 \le \cdots \le v_n \le k$.

For example, the two paths in \Cref{fig:2tuple} correspond to vectors $[1,1,2,3,6,6]$ and $[1,1,2,4,6,7]$,
which are the rows of the example matrix in \Cref{eq:matrix-example}.

Let $P_1$ denote a path from $a_1$ to $b_1$, $P_2$ a path from $a_2$ to $b_2$,
and $V_1 = [v_{1,1}, v_{1,2}, \ldots, v_{1,n}]$ and $V_2 = [v_{2,1}, v_{2,2}, \ldots, v_{2,n}]$ their corresponding row vectors.
Next, we justify why our positioning of source/destination vertices guarantee that,
if $P_1$ and $P_2$ do not intersect,
then $V_1 \le V_2$ (where $\le$ denote conjunctive component-wise $\le$),
which implies that $[V_1; V_2]$ is weakly increasing.
Consider the $i$-th upward edge in $P_1$, which should go from $(x,i-1)$ to $(x,i)$ for some $x$;
notice that $v_{1,i} = 1+x$, because the number of rightward edges in $P_1$ up to $(x,i-1)$ is exactly $x$.
Similarly, consider the $i$-th upward edge in $P_2$, which should go from $(x',i-2)$ to $(x', i-1)$ for some $x'$;
notice that $v_{2,i} = x'$, because the number of rightward edges in $P_2$ up to $(x', i-2)$ is exactly $x'-1$.
Comparing $(x,i-1)$ on $P_1$ with $(x',i-1)$ on $P_2$,
$(x,i-1)$ must be to the left of $(x',i-1)$ because $P_1$ and $P_2$ do not intersect and $P_1$ originates to the left of $P_2$.
Therefore, $x+1 \le x'$, which means that $v_{1,i} = 1+x \le x' = v_{2,i}$.

Conversely, we show that if $P_1$ and $P_2$ intersect, then $[V_1; V_2]$ is not weakly increasing.
Consider the vertex $(x,y)$ where $P_1$ and $P_2$ first intersects starting from their sources.
Note that $y \le n-1$ because $P_2$ cannot reach $y$-coordinates higher than $n-1$.
Because the paths do not intersect at any point below or to the left of this point,
$P_1$ must enter this intersection with a rightward edge $e_1$ from $(x-1,y)$ to $(x,y)$,
and $P_2$ must enter with an upward edge $e_2$ from $(x,y-1)$ to $(x,y)$.
Note the $(y+1)$-th upward edge of $P_1$ comes after $(x,y)$, so its $x$-coordinate is least $x$;
therefore, $v_{1,y+1} \ge 1+x$.
On the other hand, $e_2$ is the $(y+1)$-th upward edge of $P_2$,
so $v_{2,y+1} = x$, because the number of rightward edges in $P_2$ up to $(x,y-1)$ is exactly $x-1$.
In conclusion, $v_{1,y+1} \ge 1+x > x = v_{2,y+1}$, which means that $[V_1; V_2]$ is not weakly increasing.

Putting the above arguments together gives us the following theorem (proof is omitted because it follows from the above):
\begin{theorem}\label{thm:paths}
  There exists a bijection between $\mathfrak{M}_{n, k}$ and the set of pairs of non-intersecting lattice paths
  from $a_1 = (0,0)$ to $b_1 = (k-1,n)$ and from $a_2 = (1,-1)$ to $b_2 = (k,n-1)$.
\end{theorem}

\subsection{Counting Non-Intersecting Lattice Paths}
\label{sec:paths:lindstrom}

Given \Cref{thm:paths}, what remains to establish \Cref{eq:count-formula} is to count the number of non-intersecting lattice path pairs.
The Lindstr\"{o}m-Gessel-Viennot Lemma provides a general way to count the number of tuples of non-intersecting paths on a directed graph.
For a set of source vertices $A=\{a_1, \cdots, a_n\}$ and a set of destination vertices $B=\{b_1,\cdots,b_n\}$,
where $\varepsilon(a,b)$ is the number of paths from $a$ to $b$, the determinant of the matrix
\begin{align*}\label{eq:lindstrom}
L={\begin{bmatrix}\varepsilon(a_{1},b_{1})&\varepsilon(a_{1},b_{2})&\cdots &\varepsilon(a_{1},b_{n})\\\varepsilon(a_{2},b_{1})&\varepsilon(a_{2},b_{2})&\cdots &\varepsilon(a_{2},b_{n})\\\vdots &\vdots &\ddots &\vdots \\\varepsilon(a_{n},b_{1})&\varepsilon(a_{n},b_{2})&\cdots &\varepsilon(a_{n},b_{n})\end{bmatrix}}
\end{align*}
is equal to the number of $n$-tuples of non-intersecting paths connecting $a_i \in A$ to $b_i \in B$ for every $i \in [1,n]$.

To apply this lemma, we need to count the number of paths between any source and destination vertices in the lattice described in \Cref{sec:paths:lattice}.
Consider a source vertex $a = (x,y)$ and a destination vertex $b = (x+\Delta x, y+\Delta y)$.
A path from $a$ to $b$ consists of exactly $\Delta x$ rightward edges and $\Delta y$ upward edges;
therefore, $\varepsilon(a,b)= {\Delta x + \Delta y \choose \Delta x}$,
since we are essentially choosing $\Delta x$ rightward moves out of the $\Delta x + \Delta y$ total moves.
Thus, the matrix $L$ for our setup of \Cref{sec:paths:lattice} is:
\begin{align*}
L
={\begin{bmatrix}\varepsilon(a_{1},b_{1})&\varepsilon(a_{1},b_{2})\\\varepsilon(a_{2},b_{1})&\varepsilon(a_{2},b_{2})\end{bmatrix}}
={\begin{bmatrix}{n+k-1 \choose k-1}&{n+k-1 \choose k}\\{n+k-1 \choose k-2}&{n+k-1 \choose k-1}\end{bmatrix}}
\end{align*}
Therefore, the number of pairs of non-intersecting lattice paths is given by:
\begin{align*}
\det(L)
&= {n+k-1 \choose k-1}{n+k-1 \choose k-1}-{n+k-1 \choose k-2}{n+k-1\choose k} \\
&= {n+k-1 \choose k-1}{n+k \choose k-1}\frac{n+1}{n+k}-{n+k-1 \choose k-1}\frac{k-1}{n+1}{n+k\choose k-1}\frac{n(n+1)}{k(n+k)} \\
&= {n+k-1 \choose k-1}{n+k \choose k-1}\left[\frac{k(n+1)}{k(n+k)}-\frac{(k-1)n}{k(n+k)}\right] \\
&= \frac{{n+k-1 \choose k-1}{n+k \choose k-1}}{k}
\end{align*}
This formula matches the one for $|\mathfrak{M}_{n, k}|$ in \Cref{eq:count-formula}.

\section{Conclusion and Extension}
\label{sec:conclude}

We have presented two alternative methods for counting the number of $2 \times n$ weakly increasing matrices with entries in $[1,k]$:
first by connecting them to Kekul\'e structures of hexagonal benzenoids,
and then by connecting them to pairs of non-intersecting paths in a rectangular lattice.
We believe that both methods can be readily generalized to the problem of counting $m \times n$ weakly increasing matrices.

The first method would analyze the positions of v-bars in $\Kekule_{n,m,k-1}$.
We postulate that, going from top to bottom in the hexagonal benzenoid,
the number of v-bars in each row will first increase one at a time starting from $1$ in the top row, until it reaches $m$,
and remains at $m$ until it then decreases one at a time to $1$ in the bottom row.
The arrangement of the v-bar will follow the same staggered rightward-shifting pattern as those in $\Kekule_{n,2,k-1}$.

The second method would consider $m$ non-intersecting paths in a rectangular lattice.
Here, we further extend the set of source and destination vertices in the same way as in the case of $m=2$.
Going from one pair to the next, we shift the locations of these vertices by $(1,-1)$,
so that all source vertices lie along a line of slope $-1$ (as do all destination vertices).

\subsection*{Acknowledgments}

I would like to thank Rex Chen and Shivani Ramkumar, my fellow students at the North Carolina School of Science and Mathematics (NCSSM), for their collaboration on~\cite{ChenRamkumarYang25},
which found the explicit formula for $|\mathfrak{M}_{n, k}|$ that began my search for the two alternative approaches presented in this paper.
I would like to acknowledge Drs.\ Ashley Tharp and Jun Yang for suggesting the possibility of creating a bijection to lattice paths, which led to my solution in \Cref{sec:paths}.
I would like to especially thank Dr.\ Tharp, who was my mentor in the Mathematics Research in Mathematics Program at NCSSM,
for her extensive advice and support throughout my research process.


\newpage
\appendix
\section{Proofs}
\label{app:proof}

\begin{proof}[Proof of \Cref{lemma:top-row}]
If there are no v-bars in the top row of $K$, the edge incident to the top of unselected vertical edge $(0,0)$ must be selected, since the top vertex of $(0,0)$ must be incident to a selected edge. Along the shortest path from the top of $(0,0)$ to $(0,n)$, it follows that the edges must alternate between selected and unselected. However, this means that the top edge incident to $(0,n)$ will be unselected and the top vertex of $(0,n)$ will not be incident to a selected edge. Therefore, there must be at least one v-bar in the top row of $K$.

If some two edges $(0,a)$ and $(0,b)$ with $a < b$ in the top row are selected with all vertical edges between them being unselected, it follows that the $4$ edges that are incident to the tops of these two v-bars are unselected. From here, the edges must alternate between unselected and selected in the rightward direction from the top of $(0,a)$, and in the leftward direction from the top of $(0,b)$. However, by symmetry, at some point between $(0,a)$ and $(0,b)$, the two edges incident to this point must both be selected or unselected. If these two unselected edges are incident at this point, then there must be a v-bar at this point, which is between $(0,a)$ and $(0,b)$, creating a contradiction.

Thus, there must be one and only one v-bar in the top row of $K$. All other vertical edges in this row are unselected, and the $4$ edges incident to this v-bar are unselected. Given these restrictions, we know that along the shortest path from the top of $(0,0)$ to $(0,n)$, edges alternate between unselected and selected, starting from the $2$ unselected edges incident to the top of the v-bar, forming a willow that spans $[0,n]$.
\end{proof}

\begin{proof}[Proof of \Cref{lemma:one-to-right}]
This proof follows a logic similar to that of \cref{lemma:top-row}. Since $(i-1,j')$ is the rightmost v-bar in row $i-1$ of $K$, it follows that the two edges incident to the bottom of $(i-1,j')$ and all vertical edges to the right of it are unselected. This means that the v-bars in rows $i \in [2,r-1]$ do not determine whether the edges on the shortest path from the top of $(i,j')$ to the top of $(i,n+1)$ are selected or unselected. If $i \in [2,r-1]$, we have a row of units between $(i,j')$ and $(i,n+1)$ and at least one additional row of units below it. We cannot determine the selection of the edges in these two rows based on the composition of selected/unselected edges in rows $0$ through $i-1$. This scenario with the row between $(i,j')$ and $(i, n+1)$ exactly resembles our situation in \cref{lemma:top-row}. By \cref{lemma:top-row}, there is one and only one v-bar in this row, forming a willow anchored at $(i,j)$ and spanning $[j',n+1]$ where $j \ge j'$.
\end{proof}

\begin{proof}[Proof of \Cref{lemma:row-1}]
We can use \cref{lemma:one-to-right} to determine that if the v-bar in the top row of $K$ is in $(0,j)$, there is exactly one v-bar at $(1,j_2)$ with $j \le j_2$ that anchors a willow spanning $[j+1,n+1]$. However, unlike rows $[2,r-1]$, the row above row $1$ is the top row, which contains one fewer hexagon. If we have $j=j_2$, then the v-bar at $(1,j_2)$ would be vertically to the left of the v-bar at $(0,j)$. Thus, we would have one v-bar at $(1,j_2)$ with $j < j_2$.

Now, consider all the vertical edges to the left of $(1,j+1)$. Similarly to the situation in \cref{lemma:one-to-right}, we cannot determine whether the shortest path from the top of $(1,0)$ to the top of $(1,j)$ is selected or unselected based on the v-bar at $(0,j)$. This means that we have a row of units $(1,0)$ and $(1,j)$ and at least one row of units below it, for which we cannot determine the selection of edges (based on the top row's v-bar). This exactly resembles our situation in \cref{lemma:top-row}. By \cref{lemma:top-row}, there is one and only one v-bar between $(1,0)$ and $(1,j)$, forming a willow anchored at $(1,j_1)$ and spanning $[0,j]$ where $j_1 \le j$.
\end{proof}

\begin{proof}[Proof of \Cref{lemma:nothing-on-left}]
First, consider the base case where $i=2$. By \cref{lemma:row-1}, the leftmost v-bar in row $1$ of $K$ is $(i-1,j')$, which forms a willow spanning $[0,j'']$ where $j''>j$. If the v-bar at $(1,0)$ is selected, there cannot be a v-bar in row $2$ vertically left of $(1,0)$: we would have $j<0$. If the vertical edge $(1,0)$ is unselected, then the edge incident to the bottom of $(1,0)$ must be selected, and it follows that the shortest path from the bottom of $(1,0)$ to the top of $(2,j'-1)$ must alternate between selected and unselected edges. This means that there cannot be a v-bar between $(2,0)$ and $(2,j'-1)$. Thus, there cannot be a v-bar $(2,j)$ such that $j<j'$, and a caterpillar is formed, anchored at $(1,j')$.

Next, we have the inductive step. For some $i \in (2,r]$, assume that there is a caterpillar anchored by $(i-2,j'')$, the leftmost v-bar in row $i-2$ of $K$. Also, say the leftmost v-bar in row $i-1$ of $K$ is $(i-1,j')$. Since the vertical edge $(i-1,0)$ must be unselected (it is part of a caterpillar), the edge incident to the bottom of $(i-1,0)$ must be selected. Then it follows that the shortest path from the bottom of $(i-1,0)$ to the top of $(i,j'-1)$ must be alternating between selected and unselected edges. This means that there cannot be a v-bar $(i,j)$ such that $j<j'$, and a caterpillar is formed, anchored at $(i-1,j')$. By induction, if $(i-1,j')$ is the leftmost v-bar in row $i-1$ of $K$ for any $i \in [2,r]$, there exists no v-bar $(i,j)$ in row $i$ of $K$ such that $j < j'$, and $K$ contains a caterpillar anchored at $(i-1,j')$.
\end{proof}

\begin{proof}[Proof of \Cref{lemma:middle-rows}]
By \cref{lemma:nothing-on-left}, there will be a caterpillar in row $i$ of $K$ anchored at $(i-1,j'_1)$. By \cref{lemma:one-to-right}, there will be a willow in row $i$ of $K$ anchored at $(i,j_2)$, for some $j_2 \geq j'_2$, spanning $[j'_2,n+1]$. 

Now, note that we cannot determine whether the edges of the shortest path from the top of $(i,j'_1)$ to the top of $(i,j'_2-1)$ are selected/unselected, since they lie between (and not including) the right-hand edge incident to the bottom of v-bar $(i-1,j'_1)$ and the left-hand edge incident to the bottom of v-bar $(i-1,j'_2)$. This means that we have a row between $(i,j'_1)$ and $(i,j'_2-1)$ and at least one row below it, for which the selection of edges cannot be determined by the v-bars in rows $0$ through $i-1$. This exactly resembles our situation in \cref{lemma:top-row}, so by \cref{lemma:top-row}, there is one and only one v-bar between $(i,j'_1)$ and $(i,j'_2-1)$, which forms a willow spanning $[j'_2,n+1]$ anchored by some $(i,j_2)$ where $j_1<j'_2 \le j_2$.
\end{proof}

\begin{proof}[Proof of \Cref{lemma:bottom-row}]
Note that there must be an upside-down willow anchored at $(r,j)$. The argument is exactly the same as the proof of \Cref{lemma:top-row}, because the hexagonal benzenoid structure is symmetric by a rotation of $180$ degrees. Due to this symmetry, we can also deduce that $j'_1 \le j < j'_2$ by \cref{lemma:one-to-right}. 

Now, we need to prove that there is one possible selection of the edges in the shortest path from the bottom of $(r-1,0)$ to the bottom of $(r-1,n+1)$. First, consider the shortest path from the bottom of $(r-1,0)$ to the bottom of $(r-1,j'_1)$. There will always be an even number of edges in this path (including $0$). The edge incident to the bottom of $(r-1,0)$ must be selected and the edge incident to the bottom of $(r-1,j'_1)$ must be unselected, so there is only one possible selection of edges in this path (if $j'_1=0$, there is also only one arrangement: there are no edges to select). Following a similar proof, there is also only one selection of edges in the shortest path from the bottom of $(r-1,j'_2)$ to the bottom of $(r-1,n+1)$.

Next, consider the shortest path from the bottom of $(r-1,j'_1)$ to the top of $(r,j)$. There will always be an odd number of edges in this path. The edge incident to the bottom of $(r-1,j'_1)$ and the edge incident to the top of $(r,j)$ must always be unselected, so there is only one possible selection of edges in this path. Following a similar proof, there is also only one selection of edges in the shortest path from the top of $(r,j)$ to the bottom of $(r-1,j'_2)$. Thus, the shortest path from the bottom of $(r-1,0)$ to the bottom of $(r-1,n+1)$ must be selected/unselected in a unique way.
\end{proof}

\begin{proof}[Proof of \Cref{lemma:benzenoid}]
The proof follows directly from lemmas \ref{lemma:top-row} through \ref{lemma:bottom-row}.
\end{proof}

\begin{proof}[Proof of \Cref{lemma:pulse}]
Every $M \in \mathfrak{M}_{n,k}$ can be expressed as the sum of a $2 \times n$ matrix with every entry being a $1$ and a unique combination of $k-1$ pulse matrices of $n$ columns, in decreasing order. These pulse matrices can be determined first by subtracting $J_{2,n}$ from $M$. Then, look at both rows of this resultant matrix. To obtain a pulse matrix of $n$ columns based on this resultant matrix, assign all nonzero entries in the resultant matrix to $1$ and all zero entries to $0$. Take note of this pulse matrix (this will be the largest pulse matrix) and subtract it from the resultant matrix: now consider this new resultant matrix. Repeat this process, obtaining a new pulse matrix from the nonzero entries of the resultant matrix and subtracting it from the resultant matrix, until we have extracted $k-1$ pulse matrices and the resultant matrix is the zero matrix. It follows that the $k-1$ pulse matrices obtained will be weakly increasing along the rows and columns and be generated in order of decreasing magnitude (Frobenius norm). Thus, every $M \in \mathfrak{M}_{n,k}$ can be written uniquely as $J_{2,n} + \sum_{i=1}^{k-1} L_i$, where each $L_i$ is a pulse matrix, and $L_{i-1} \ge L_i$ for all $i \in [2,k-1]$.
\end{proof}

\begin{proof}[Proof of \cref{thm:benzenoid}]
By lemmas \ref{lemma:top-row} through \ref{lemma:benzenoid}, each Kekul\'e structure in $\Kekule_{n,2,k-1}$ can be divided into $k-1$ groups of two adjacent rows, with each group containing $2$ v-bars. The first group of two rows contains the v-bar in the $0$-th row and the leftmost v-bar in the $1$st row, and every successive group of two rows $r$ and $r+1$ contains the rightmost v-bar of row $r$ and the leftmost v-bar of row $r+1$. By letting the v-bars of these groups of $2$ adjacent rows denote the locations of the two leftmost $1$s in a pulse matrix, a Kekul\'e structure can be expressed as a unique combination of $k-1$ pulse matrices in order of decreasing magnitude (Frobenius Norm) (\cref{lemma:benzenoid}). Summing these pulse matrices and $J_{2,n}$ yields a unique weakly increasing matrix in $\mathfrak{M}_{n,k}$.

By \cref{lemma:pulse}, every matrix $M \in \mathfrak{M}_{n,k}$ can be written as a unique sum of $J_{2,n}$ and $k-1$ pulse matrices in order of decreasing magnitude. Letting the positions of the leftmost $1$s of the first and second rows of these matrices denote the positions of the $2$ v-bars in each group of two adjacent rows, we can express every matrix $M \in \mathfrak{M}_{n,k}$ as a unique Kekul\'e structure in $\Kekule_{n,2,k-1}$. Thus, we have drawn a bijection between every matrix $M \in \mathfrak{M}_{n,k}$ and every Kekule structure $K \in \Kekule_{n,2,k-1}$.
\end{proof}

\end{document}